\documentclass{amsart}
\usepackage{amssymb}
\usepackage[hyphens]{url} 

\theoremstyle{plain}

\theoremstyle{definition}

\theoremstyle{remark}
 
 \numberwithin{equation}{section}

\renewcommand{\leq}{\leqslant}

\renewcommand{\setminus}{\smallsetminus}
\usepackage{amsmath}
\newcommand{\legendre}[2]{\left(\frac{#1}{#2}\right)}

\title[On distinct residues of factorials]{ON DISTINCT RESIDUES OF FACTORIALS}



\author[Andreji\'c]{\bfseries Vladica Andreji\'c}

\address{
Faculty of Mathematics \\
University of Belgrade  \\
Belgrade\\
Serbia} 
\email{andrew@matf.bg.ac.rs}

\author[Tatarevic]{\bfseries Milos Tatarevic}
\address{Alameda, CA 94501}
\email{milos.tatarevic@gmail.com}

\thanks{This work is partially supported by the Serbian Ministry of Education and Science, project No. 174012}

\begin{document}

\setcounter{page}{1}
\thispagestyle{empty}

\begin{abstract}
We investigate the existence of primes $p > 5$ for which the residues of
$2!$, $3!$, \dots, $(p-1)!$ modulo $p$ are all distinct. We describe the connection
between this problem and Kurepa's left factorial function, and report that there
are no such primes less than $10^{11}$.
\end{abstract}

\maketitle

\section{Introduction}

Paul Erd\H{o}s asked the following question: Is there a prime $p > 5$ such that the residues
of $2!$, $3!$, \dots, $(p-1)!$ modulo $p$ are all distinct?
Rokowska and Schinzel \cite{RS} proved that if such $p$ exists, it needs to satisfy
the following conditions

\begin{equation}\label{RS-c}
p\equiv 5\;(\bmod\; 8),\quad
\legendre{5}{p}=-1,\quad
\legendre{-23}{p}=1,
\end{equation}

\noindent and that the missing residue must be that of $-((p-1)/2)!$.
The given conditions enabled them
to prove that there are no such primes $p$ with $5<p<1000$.
This problem is also mentioned in \cite[Section F11]{Gu}.
Recently, Trudgian \cite{T} called such a prime $p$ a socialist
prime and proved that $p$ also needs to satisfy

\begin{equation}\label{T-c}
\begin{aligned}
\legendre{1957}{p}=1, & \text{ or } \legendre{1957}{p}=-1 \text{ with } \legendre{4y+25}{p}=-1 \\
 & \text{ for all $y$ satisfying } y(y+4)(y+6)-1\equiv 0\pmod p.
\end{aligned}
\end{equation}
He confirmed that there are no such primes less than $10^9$.

In this paper, we describe the connection between the socialist primes and the left factorial
function $!n=0!+1!+\cdots+(n-1)!$ introduced by \DJ uro Kurepa \cite{K1}. Kurepa
conjectured that $\text{GCD}(!n,n!)=2$ holds for all integers
$n>1$, which is equivalent to the statement that there is no odd
prime $p$ that divides $!p$. This conjecture is also
mentioned in \cite[Section B44]{Gu}.

In our previous work \cite{AT}, we calculated and recorded the residues $r_{p}\!=\,\,!p \bmod p$ for all
primes $p<2^{34}$. Now we show that if $p$ is a
socialist prime then $(!p-2)^2\equiv -1 \pmod p$, which enabled us
to immediately confirm that there are no such primes
less than $2^{34}$. Additionally, we extended the search up to $10^{11}$.

\section{Left factorial calculations}

Let us demonstrate some straightforward calculations. Wilson's
theorem states that
\begin{equation}\label{wil1}
(p-1)!\equiv -1\pmod p,
\end{equation}
for all primes $p$. Therefore,
\begin{equation}\label{wil2}
(p-2)!\equiv 1\pmod p,
\end{equation}
and more generally
\begin{equation}\label{wilson}
(p-k)!(k-1)!\equiv(-1)^k \pmod p,
\end{equation}
for all primes $p$ and all $1\leq k\leq p$.
Especially, we have $(((p-1)/2)!)^2\equiv (-1)^{(p+1)/2} \pmod p$.
Since in (\ref{wil1}) and (\ref{wil2}) we already have residues $-1$ and $1$, if $p$ is a socialist prime it follows that
$((p-1)/2)!\not\equiv \pm 1 \pmod p$, so
$(-1)^{(p+1)/2}\neq 1$, thus $p\equiv 1 \pmod 4$, and consequently
\begin{equation}\label{wil3}
\left(\left(\frac{p-1}{2}\right)!\right)^2\equiv -1 \pmod p.
\end{equation}

Let $r$ be the missing residue. Then we have
\begin{equation*}
(p-1)!\equiv r \prod_{k=2}^{p-1}k!=r(p-2)!(p-1)!\left(\frac{p-1}{2}\right)! \prod_{k=3}^{(p-1)/2}(p-k)!(k-1)! \pmod p,
\end{equation*}
which after (\ref{wil1}), (\ref{wil2}), (\ref{wilson}), and (\ref{wil3}) leads to
$r\equiv (-1)^{(p^2-1)/8} ((p-1)/2)! \pmod p$.
Since $r\not\equiv ((p-1)/2)! \pmod p$, we have
\begin{equation}\label{res-r}
r \equiv-\left(\frac{p-1}{2}\right)! \pmod p.
\end{equation}
Consequently, $(p^2-1)/8$ is odd, and therefore, $p\equiv 5 \pmod 8$.

Since (\ref{res-r}), we can connect a socialist prime $p$ with the left factorial function. We have
\begin{equation*}
-\left(\frac{p-1}{2}\right)!+\sum_{k=2}^{p-1}k!\equiv \sum_{i=1}^{p-1}i=\frac{(p-1)p}{2}\equiv 0 \pmod p,
\end{equation*}
which implies
\begin{equation*}
\left(\frac{p-1}{2}\right)!\equiv \sum_{k=2}^{p-1}k!=\sum_{k=0}^{p-1}k!-2= !p -2 \pmod p.
\end{equation*}
Finally (\ref{wil3}) gives the necessary condition

\begin{equation}\label{lfc}
(!p -2)^2\equiv -1 \pmod p.
\end{equation}

In our previous work \cite{AT}, we calculated and recorded the residues
$r_{p}\!=\,\,!p\pmod p$ for all primes $p<2^{34}$.
After a fast search through our database, the only primes
where $p$ divides $(r_p-2)^2+1$ with $p<2^{34}$ are $5$, $13$,
$157$, $317$, $5449$, and $5749$. More concretely, $r_5=4$,
$r_{13}=10$, $r_{157}=131$, $r_{317}=205$, $r_{5449}=4816$, and
$r_{5749}=808$. Consequently, there are no socialist primes $p$ with
$5<p<2^{34}$. It is interesting that there are no small $p\equiv
3\pmod 4$ with $p \mid (r_p-2)^2+1$.

In \cite{AT}, we also considered the generalized Kurepa's left factorial
\begin{equation*}
!^k n=(0!)^k+ (1!)^k+\cdots+((n-1)!)^k,
\end{equation*}
where $!^1 k=!k$. There we presented counterexamples for the analog of Kurepa's conjecture, 
where for all $1<k<100$ there exists an odd prime $p$ such that $p\mid\,!^k p$.  
Similarly, we can connect the socialist primes with a generalized left factorial in the following way. 
From
\begin{equation*}
p^{n+1}-1= \sum_{m=1}^{p-1}\left((m+1)^{n+1}-m^{n+1}\right)=
\sum_{m=1}^{p-1}\sum_{k=0}^n \binom{n+1}{k}m^k,
\end{equation*}
we have
\begin{equation*}
p^{n+1}-1=p-1+\sum_{k=1}^n\binom{n+1}{k}\sum_{m=1}^{p-1}m^k,
\end{equation*}
and therefore,
\begin{equation*}
\sum_{k=1}^n \binom{n+1}{k}(1^k+2^k+...+(p-1)^k)\equiv 0\pmod p
\end{equation*}
for all $n$,
which gives $1^k+2^k+...+(p-1)^k \equiv 0\pmod p$ for all $1\leq
k\leq p-2$. Thus $$(2!)^k+...+((p-1)!)^k
+\left(-\left(\frac{p-1}{2}\right)!\right)^k\equiv
1^k+2^k+...+(p-1)^k \equiv 0\pmod p$$ and hence
\begin{equation*}
!^k p=(0!)^k+...+((p-1)!)^k\equiv 2-\left(-\left(\frac{p-1}{2}\right)!\right)^k \pmod p.
\end{equation*}
If we include (\ref{wil3}), we finally get the necessary condition:

\begin{equation}\label{lfck}
\begin{split}
(!^k p-2)^2+1\equiv 0 \pmod p &\quad \mbox{if $k$ is odd},\\
!^k p\equiv 1 \pmod p &\quad \mbox{if $k=4t$},\\
!^k p\equiv 3 \pmod p &\quad \mbox{if $k=4t+2$}.
\end{split}
\end{equation}

\section{Additional calculations}

Let us consider the set $H=\{2,3,4,...,p-3\}\setminus\{\frac{p-1}{2}\}$ and let $p$ be a socialist prime. 
There is a function $f$ defined on $H$, such that for all $k\in H$,
\begin{equation}\label{ef}
(f(k))!\equiv -k!\pmod p.
\end{equation}
Thus, $f$ is an involution on the set $H$. By using (\ref{wilson})
after we multiply (\ref{ef}) by $(p-1-f(k))!(p-1-k)!$ we
get
\begin{equation}
(-1)^{f(k)+1}(p-1-k)!\equiv -(-1)^{k+1}(p-1-f(k))!\pmod p.
\end{equation}
Since $(p-1-k)!\not\equiv (p-1-f(k))! \pmod p$, we conclude that
$f(k)\equiv k\pmod 2$. This splits the set $H$ into
$(p-5)/4$ quadruples with a shape
$$U_{k}=\{k,f(k),p-1-k,p-1-f(k)\}$$
with $\prod_{x\in U_k}x\equiv 1$,
$\sum_{x\in U_k}x\equiv 0\pmod p$, and $x\equiv y\pmod 2$ for all
$x,y\in U_k$.

One idea can be related to $!^k p$ for not fixed $k$, for example $k=(p-1)/2$.
Since
\begin{equation*}
\legendre{x!}{p}=\legendre{-(f(x)!)}{p}=\legendre{-1}{p}\legendre{f(x)!}{p}=\legendre{f(x)!}{p}, 
\legendre{x!}{p}=\legendre{\frac{1}{x!}}{p}, 
\end{equation*}
we see that all members of $U_{k}$ have the same quadratic residue modulo $p$.
In this quadruple two members are less than $(p-1)/2$, and therefore,
\begin{equation*}
\legendre{2!\cdot 3! \cdots \frac{p-3}{2}!}{p}=1,
\end{equation*}
or more strict $\legendre{2!\cdot 4! \cdots ((p-5)/2)!}{p}=1=\legendre{3!\cdot 5! \cdots ((p-3)/2)!}{p}$.
Then we have 
\begin{equation*}
1=\legendre{2!\cdot 3! \cdots \frac{p-3}{2}!}{p}=\legendre{3\cdot 5 \cdots \frac{p-3}{2}}{p}=
\legendre{\frac{p-3}{2}!!}{p}= 
\legendre{({\frac{p-1}{2}!})/({2^{\frac{p-1}{4}} \frac{p-1}{4}!})}{p}.
\end{equation*}
Since $\legendre{((p-1)/2)!}{p}=(((p-1)/2)!)^{(p-1)/2}=(-1)^{(p-1)/4}=-1$ and
$\legendre{2^{(p-1)/4}}{p}=\legendre{2}{p}^{(p-1)/4}=\legendre{2}{p}=(-1)^{({p^2-1})/8}=-1$,
we can conclude that
\begin{equation*}
\legendre{\frac{p-1}{4}!}{p}=1.
\end{equation*}

\section{Heuristic Considerations}

Let us suppose that factorials $2!$, $3!$, ..., $(p-1)!$ modulo $p$ are random nonzero integers.
The probability that $p$ is a socialist prime can be estimated by
\begin{equation*}
1\times \frac{p-2}{p-1}\times \frac{p-3}{p-1}\times\dots\times\frac{2}{p-1}=\frac{(p-2)!}{(p-1)^{p-3}}.
\end{equation*}
If we include the fact that $(k+1)!\not\equiv k!$ for $2\leq k\leq p-2$ then the estimated probability is slightly higher
\begin{equation*}
W_p=1\times \frac{p-2}{p-2}\times \frac{p-3}{p-2}\times\dots\times\frac{2}{p-2}=\frac{(p-2)!}{(p-2)^{p-3}}.
\end{equation*}
From Stirling's approximation, we have $k! \approx \sqrt{2\pi k} \binom{e}{k}^k \leq e k^{k+\frac{1}{2}}e^{-k}$, and consequently
\begin{equation}\label{wp}
W_p\leq (p-2)^{\frac{3}{2}}e^{3-p}.
\end{equation}

Note that this is just a rough upper bound. If the conditions (\ref{RS-c}) and (\ref{T-c}) are included,
this estimation can be reduced by some factor. If we make an assumption that any $!^k p \pmod p$
for $k=1,\dots,p-2$ is an independent and random number, by including the condition (\ref{lfck})
one can conclude that $W_p \approx p^{2-p}$. 
However, we have that $!^{2k} p \equiv\,\, !^{p-2k-1} p \pmod p$ for $1\leq k\leq (p-3)/2$ and odd primes $p$,
which implies that $W_p$ should be greater than previously assumed.

Further, we can estimate the number of socialist primes in an interval $[a,b]$ as a
sum of $W_p$ over the primes
\begin{equation*}
\sum_{a\leq p \leq b}W_p\approx \sum_{n=a}^{b}W_{n\ln n} \approx \int_{a}^{b} W_{t\ln t} \,dt.
\end{equation*}
From (\ref{wp}), we can see that $W_{p}< e^3 p^{\frac{3}{2}}e^{-p}$, therefore,
\begin{equation*}
\begin{split}
\sum_{a\leq p \leq b}W_p & <  e^3 \int_{a}^{b} t^{\frac{3}{2}-t}(\ln t)^{\frac{3}{2}} \,dt\\
& < e^3 \int_{a}^{b} t^{\frac{3}{2}-t}(\ln t)^{\frac{3}{2}} \left(1+\frac{1}{\ln t}-\frac{3}{2t\ln t}-\frac{1}{2t(\ln t)^2} \right)\,dt\\
& < e^3 \int_{a}^{b} \left(-t^{\frac{3}{2}-t}(\ln t)^{\frac{1}{2}}  \right)'\,dt
\end{split}
\end{equation*}
and thus
\begin{equation*}
\sum_{a\leq p \leq b}W_p< e^3a^{\frac{3}{2}-a}\sqrt{\ln a}.
\end{equation*}

According to the last statement, we expect no more than $e^3a^{\frac{3}{2}-a}\sqrt{\ln a}$ socialist
primes greater than $a$.
As we confirmed that there are no socialist primes less than $10^{11}$, the estimated probability that 
such primes exist is less than $10^{-10^{12}}$.

\section{Computer search}

In 1960, Rokowska and Schinzel \cite{RS} reported that there are no primes $p$ with $5<p<1000$
for which the residues of $2!$, $3!$, \dots, $(p-1)!$ modulo $p$ are all distinct.
By applying the conditions (\ref{RS-c}), there are only ten primes below 1000 that need to be examined. 
For $5<p<10^{6}$, there are only $4908$ such primes, and after applying (\ref{T-c}) this is further
reduced to $3662$ primes \cite{T}. Recently, Trudgian \cite{T} confirmed that there are no such primes
less than $10^9$.

As a part of the search for a counterexample to Kurepa's conjecture, we recorded the residues
$r_{p}\!=\,\,!p \bmod p$ for all $p<2^{34}$ \cite{AT}. By using the congruence (\ref{lfc}), we
instantly verified there are no such primes $p$ that satisfy this condition for $10^{9} < p < 2^{34}$.
To extend the range beyond $2^{34}$ we need a more efficient method, as the time complexity of
the algorithm to obtain $r_{p}$ for all $p < n$ is $O(n^{2}/\ln n)$.

To show that $p$ is not a socialist prime, it is sufficient to find a single pair of integers $i$ and $j$,
such that $2\leq i < j \leq p - 1$ and $i! \equiv j! \pmod p$.
As we can consider that $i!$ and $j!$ modulo $p$ are ``random" integers, to find such a pair we
can apply a well-known probabilistic method called the birthday attack.
On average, it is expected that the first such a pair will be found after approximately $\sqrt{p\pi / 2}$
attempts \cite{KN}. The time complexity of this algorithm is $O(\sqrt{p})$ for single $p$,
and $O(n^{{3}/{2}}/\ln(n))$ for all $p < n$.
As we already mentioned, we do not need to examine all the primes.
To simplify the calculation, we used congruences given in (\ref{RS-c}). Note that in
this case, the time complexity remains the same as the search space is reduced by the constant
factor. Using this method, we confirmed there are no socialist primes less than $10^{11}$.

For our computation, we used a single Intel Core i7-4980HQ CPU running at 2.8GHz.
This CPU natively supports $64$-bit multiplication without overflow, which allows us to implement
fast modular reduction.
Note that we do not necessarily need to search for the first occurrence of two the same factorials
modulo $p$. In our implementation, we made several changes that, on average, slightly increased
the number of iterations required to register these duplicates, but overall they made the execution faster.

To reduce the memory consumption, we used an array of integers as a simple hash table
without collision resolution.
Accessing the elements of the hash table can be slow if accessed memory blocks are not in the cache.
The CPU we used has $256$KB of L$2$ cache per core and $6$MB of shared L$3$ cache.
While we focus on primes less than $10^{11}$, we can rarely expect more
than $4 \cdot 10^5$ iterations before the first duplicate is discovered. This number of iterations
is large enough to assume that it will be hard to keep all the data in the L$2$ cache.
Using the right size of the array we targeted the efficient usage of
L$3$ cache. Various sizes of arrays were tested, and the best parameters were determined
empirically. In particular, for $p$ in the region of $10^{11}$, an array of $2^{19}$ elements
provided good results.

The search for all $p < 10^{11}$ took slightly over one day. Although the search can be
extended beyond this bound, it is reasonable to believe that there are no socialist primes, as we
explained in Section 4.

\bibliographystyle{amsplain}

\end{document}